\documentclass[10pt]{article}
\usepackage{amssymb, amsthm}

\newcommand{\Z}{\mathbb{Z}}

\newcommand{\x}{\ma{x}}

\newcommand{\e}{\emph}
\newcommand{\rom}{\mathrm}
\newcommand{\bfP}{\mathbb{P}}

\newcommand{\ma}{\mathbf}
\newcommand{\ben}{\begin{enumerate}}
\newcommand{\een}{\end{enumerate}}
\newcommand{\eit}{\begin{itemize}}
\newcommand{\beq}{\begin{equation}}
\newcommand{\eeq}{\end{equation}}

\newcommand{\ve}{\varepsilon}

\newcommand{\mcal}{\mathcal}
\newcommand{\lb}{\left(}
\newcommand{\rb}{\right)}
\newcommand{\lab}{\label}

\newcommand{\al}{\alpha}
\newcommand{\D}{\Delta}

\newcommand{\be}{\beta}

\newcommand{\sfl}{\mathsf{\Lambda}}

\newcommand{\no}{\noindent}
\newcommand{\h}{\rom{h.c.f.}}
\newtheorem{thm}{Theorem}
\newtheorem{lem}{Lemma}

\begin{document}

\title{Plane curves in boxes and equal sums of two powers}

\author{T.D. Browning$^1$ \\ D.R. Heath-Brown$^2$\\
\small\emph{Mathematical Institute,
24--29 St. Giles',
Oxford OX1 3LB}\\
\small{$^1$browning@maths.ox.ac.uk}, \small{$^2$rhb@maths.ox.ac.uk}}
\date{}

\maketitle

\section{Introduction}

Let $F \in \Z[x_1,x_2,x_3]$ be an absolutely irreducible form of
degree $d$, producing a plane curve in $\bfP^2$.
The central aim of this paper is to analyze
the density of rational points on such curves, which are
contained in boxes with unequal sides.
We shall see below how such considerations may be used to obtain new
paucity
results for equal sums of two powers.

Suppose that $\ma{P}=(P_1,P_2,P_3)$ for fixed real numbers
$1 \leq P_1 \leq  P_2 \leq  P_3$, say.  Then we define
$$
N(F;\ma{P})= \#\{\x \in \Z^3: ~F(\x)=0, ~|x_i| \leq
P_i, ~(1\leq i \leq 3), ~\mbox{$\x$ primitive}\},
$$
where $\x=(x_1,x_2,x_3)$ is said to be primitive if
$\h(x_1,x_2,x_3)=1$.
Our starting point is the recent work of the second author
\cite{annal}, who has shown that
$N(F;\ma{P})=O(P^{2/d+\ve})$
for any $\ve>0$, whenever $P_i=P$ for $1\leq i \leq 3$.
The implied constant in this bound depends at most
upon the choice of $\ve$ and $d$, a convention that we shall follow
throughout this paper.
This is essentially best possible for curves of genus zero,
and it is natural to ask what can be said about the quantity $N(F;\ma{P})$
when the $P_i$ are of genuinely different sizes.
With this in mind we define
\beq\lab{TT}
T=\max\{P_1^{e_1}P_2^{e_2}P_3^{e_3}\},
\eeq
where the maximum is taken over all integer triples
$(e_1,e_2,e_3)$ for which the corresponding monomial
$x_{1}^{e_{1}}x_2^{e_2}x_{3}^{e_{3}}$
occurs in $F(\ma{x})$ with non-zero coefficient.
Then for any  $\ve>0$, the second author's principal result for curves
\cite[Theorem
3]{annal} states that
\beq\lab{annals}
N(F;\ma{P}) \ll \lb\frac{P_1P_2P_3}{T^{1/d}} \rb^{1/d+\ve}.
\eeq
This clearly reduces to the previous bound whenever $P_i=P$ for $1\leq i
\leq 3$.
Turning to the case of unequal $P_i$, it is easy
to construct examples for which (\ref{annals}) is not best possible.
Indeed, let $F= x_1^{d-1}x_3-x_2^d$ and $\ma{P}=(1,P,Q)$, say.  Then it
follows from (\ref{annals}) that $N(F;\ma{P})\ll
P^{1/d}Q^{(d-1)/d^2+\ve}$ whenever $P^d \leq Q$, whereas in fact
$N(F;\ma{P})$ has order of magnitude $\min\{P,Q^{1/d}\}$.

It transpires that in the case of unequal $P_i$ there is
scope for improvement within the proof of (\ref{annals}) itself.  This
has been demonstrated by the second author \cite[Theorem 15]{cime}
in the special case $P_1=1$.
For any $\ve>0$, it is shown that
\beq\lab{cime}
N(F;1,P_2,P_3) \ll  P_3^\ve
\exp\lb \frac{\log P_2\log P_3}{\log T}\rb,
\eeq
where $T$ is given by (\ref{TT}).  In particular, this is always at
least as sharp as  (\ref{annals}) and is
essentially best possible by our example above.
Moreover, the formulation we have given obviously incorporates the
corresponding problem for integral points on  affine plane curves
$F(1,x_2,x_3)=0$.

We build upon these results by returning once again to the framework
provided by the proof
of (\ref{annals}).  Our aim is to establish
a sharper bound for the interim case in which $P_1\geq 1$ and
the region $|x_i| \leq P_i$ is sufficiently lopsided.
Unfortunately the statement of the bound is somewhat
complicated, and it will be convenient to introduce the quantities
$$
\al=\frac{\log P_1}{\log P_3}, \qquad \be=\frac{\log P_2}{\log P_3},
\qquad \tau=\frac{\log T}{d\log P_3},
$$
for given $1 \leq P_1 \leq  P_2 \leq  P_3$ and $T$ as in
(\ref{TT}).
With this in mind, we have the following result.

\begin{thm}
Suppose that $T \geq (P_1P_2)^d$, and let $\ve>0$ be given.  Then we have
$$
N(F;\ma{P}) \ll  P_3^\ve\exp\lb
\frac{\al \be +(\al+ \be
-\al\be)(\tau-\al-\be)}{d(\tau-\al)(\tau-\be)}\log P_3
\rb.
$$
\end{thm}

In particular, Theorem 1 reduces to (\ref{cime}) in the case $P_1=1$.
Indeed, we
always have  the lower bound $T \geq P_2^d$ whenever $F$ is an
absolutely irreducible form, and the bounds in Theorem 1 and (\ref{cime})
agree when $\al=0.$
On the face of it, one might think of the condition $T \geq
(P_1P_2)^d$ as being unduly restrictive.  In fact a straightforward
calculation
shows that this is precisely the region for which Theorem 1 is sharper
than (\ref{annals}).

Turning to our application of Theorem 1, we fix a choice of $k
\geq 4$ and consider the diagonal
equation
\begin{equation}\lab{x}
w^{k}+x^{k}=y^{k}+z^{k}.
\end{equation}
For any $X \geq 1$,
we denote by $N_{k}(X)$ the
number of positive integer solutions
in the region $\max \{w,x,y,z\}\leq X$.
There are $2X^{2}+O(X)$
trivial solutions in which $y,z$
are a permutation of $w,x$, and so we write $N_{k}^{(0)}(X)$ for
the number of non-trivial solutions.  This quantity has received a
great deal of attention lately, and we mention in particular the
results of Hooley \cite{hoo1, hoo2} and the second author
\cite[Theorem 11]{annal}.
Together, they  comprise the best available
estimates for values of $k$ in the interval $4 \leq k \leq 12$.
The first of these provides the bound
\beq\lab{hooley}
N_{k}^{(0)}(X)\ll X^{5/3+\ve}
\eeq
for any $k \geq 4$, whereas the second yields
\beq\lab{11}
N_{k}^{(0)}(X) \ll X^{1+\ve}+X^{3/\sqrt{k}+2/(k-1)+\ve}
\eeq
for any such $k$, and supersedes Hooley's bound  for $k \geq 6$.

The aim of the second part of this paper is to improve upon the
previous bounds whenever $k=5$ or $6$.  This will be done via a suitable
application of Theorem 1.   It is unfortunate that  we shall only make
use of the special case (\ref{cime}), and not of Theorem 1 \e{per
se}.  Nonetheless, it is our belief that the bound in Theorem 1
still merits a full presentation.
We shall establish the following result in Section $3$.

\begin{thm}
For any $k\geq 4$ and any $\ve>0$, we have
$$
N_{k}^{(0)}(X) \ll X^{3/2+1/(2k-2)+\ve}.
$$
\end{thm}

We take this opportunity to remark that the proof of (\ref{11}) may be
modified slightly to give a sharper result.
In fact it is possbile to establish the estimate
\beq\lab{11'}
N_{k}^{(0)}(X) \ll X^{1+\ve}+X^{3/\sqrt{k}+2/k+\ve},
\eeq
for any $\ve>0$ and $k\geq 4$.  At this point it is convenient to
tabulate the various available bounds for $N_{k}^{(0)}(X)$, for
$k$ in the range $4 \leq k \leq 8$.
Let $\ve>0$.  Then we may write
$N_{k}^{(0)}(X)=O(X^{\theta_k+\ve})$, where the permissible values of
$\theta_k$
are given in the following table. The rows in this table
correspond to the estimates
(\ref{hooley}), (\ref{11}), (\ref{11'}) and Theorem
$2$, respectively.

\begin{table}[!h]
\begin{center}
\begin{tabular}{|l| l| l|l|l|}
\hline
$\theta_4$ & $\theta_5$ & $\theta_6$ & $\theta_7$ & $\theta_8$\\
\hline
\hline
$1.666..$ & $1.666..$ & $1.666..$& $1.666..$ & $1.666..$\\
$2.166..$ & $1.841..$ & $1.624..$ & $1.467..$ & $1.346..$\\
$2.000..$ & $1.741..$ & $1.558..$ & $1.419..$ & $1.310..$\\
$1.666..$ & $1.625$ & $1.600..$ & $1.583..$ & $1.571..$\\
\hline
\end{tabular}
\end{center}
\end{table}

\no
Thus we see that Hooley's bound (\ref{hooley}) remains unbeaten only
for
$k=4$.  For $k=5$ the exponent in Theorem $2$ is the sharpest known,
but (\ref{11'}) should be used for larger values of $k$.

We now indicate how (\ref{11'}) can be established.
An inspection
of the proof \cite[\S 8]{annal} of (\ref{11}), reveals that it
suffices to offer an alternative treatment of the curves of degree $k-1$
which are contained in the non-singular projective surface (\ref{x}).
Our observed improvement rests upon a reformulation of Colliot-Th\'el\`ene's
result \cite[Appendix]{annal}, as used in the proof of (\ref{11}).  This
states  that any
non-singular surface of degree $k$ in $\bfP^3$  contains
$O(1)$ curves of degree $\leq k-2$.  In recent communications with the
authors,
Professor Colliot-Th\'el\`ene  has
shown that any such non-singular surface actually contains
$O(1)$ non-degenerate curves of degree $\leq 2(k-2)$.  Here, a
curve in $\bfP^3$ is said to be non-degenerate if it is not contained
in any $\bfP^2 \subset \bfP^3$.  Since any absolutely irreducible curve of
degree $d$ in $\bfP^3$
contains $O(X^{2/d+\ve})$ rational points of height at most
$X$, by \cite[Theorem 5]{annal},
it remains to handle the plane curves of degree $k-1$ which are
contained in (\ref{x}).  Such curves arise as the intersection of
(\ref{x}) with a plane
\beq\lab{plane}
a w+b x+c y+d z=0,
\eeq
say, that contains one of the lines in (\ref{x}).
But we know that all of the lines contained in this surface
are given by $\{|w|,|x|\}= \{|y|,|z|\}$ (for $k$ even), or by
$\{w,x\}= \{y,z\}$ or $w=-x$, $y=-z$ (for $k$ odd). Hence it is trivial to
deduce from (\ref{plane}) that the only available planes have 
$\{|a|,|b|\}= \{|c|,|d|\}$ (for $k$ even), or $\{a,b\}= \{-c,-d\}$
or  $a=b$, $c=d$ (for $k$ odd).
Thus there are relatively few planes (\ref{plane}) of low height that
need to be considered.  The proof may then be completed by counting
points according to the height of the corresponding plane (\ref{plane}).

In the case $k=4$, it is worthwhile remarking that the proof of
Theorem $2$ can readily be adapted to show that for any $\ve>0$ there are
$O(X^{5/3+\ve})$ positive integer solutions to the equation
$$
w^4+x^4+y^4=z^4,
$$
in the region $\max \{w,x,y,z\}\leq X$.  This supersedes work of the
first author \cite{browning.spm}, who has already obtained the
exponent $7/4+\ve$.

{\bf Notation.}
We shall follow common practice in allowing
the small positive quantity $\ve$ to take different values at different
points in
all that follows.

{\bf Acknowledgement.} While working on this paper, the first
author was supported by
EPSRC Grant number GR/R93155/01.

\section{Proof of Theorem 1}

In this section we shall prove Theorem 1.
If $P_1=1$, then Theorem 1 reduces to (\ref{cime}).
Henceforth we assume that $P_1>1$.
But then the
condition $T \geq (P_1P_2)^d$ implies that $f_3
\neq 0$, where we suppose that $(f_1,f_2,f_3)$ is the maximal
triple taken in the definition (\ref{TT}) of $T$.  We set
\beq\lab{kappa}
\kappa = (2f_1+f_2)/f_3+3,
\eeq
and note that $3 \leq \kappa \leq 3d$.
During the course of our argument we will encounter difficulties if the
values
of $\log P_i$ are too close together.   We therefore replace
$P_1,P_2,P_3$ by
\beq\lab{b}
B_1=P_1P_3^{\delta}, \qquad B_2=P_2P_3^{2\delta}, \qquad
B_3=P_3^{1+\kappa\delta}.
\eeq
Here $\kappa$ is given by (\ref{kappa}), and $\delta$ is defined to be
\beq\lab{del}
\delta=\frac{\ve}{180d^3}.
\eeq
Writing  $B_1^{f_1}B_2^{f_2}B_3^{f_3}=T'$, say, we observe that
$$
\log T' =  \log T+3 d \delta \log P_3.
$$
For any $x \geq 0$, define the functions
\beq\lab{abt1}
\al(x)=\frac{\al+x}{1+\kappa x}, \qquad
\be(x)=\frac{\be+2x}{1+\kappa x}, \qquad \tau(x)=\frac{\tau+3x}{1+\kappa
x},
\eeq
where $\al, \be, \tau$ are the quantities appearing in the statement
of Theorem 1.  Then
(\ref{b})  implies that
$$
\al(\delta)=\frac{\log B_1}{\log B_3}, \qquad \be(\delta)=\frac{\log
B_2}{\log B_3},
\qquad \tau(\delta)=\frac{\log T'}{d\log B_3}.
$$
It will be convenient to record that for any $x \geq 0$ we have
$0<\al(x) \leq \be(x)$ and
\beq\lab{tau}
\al(x)+\be(x) \leq \tau(x) \leq 1,
\eeq
since $P_3^d \geq T \geq (P_1P_2)^d$ and $\kappa \geq 3$.
Finally, we define the function
\beq
g(x) = \frac{\al(x)^2}{\tau(x)-\be(x)}\cdot\frac{\be(x)
  -1}{\tau(x)-\al(x)} + \frac{\al(x)+ \be(x)
-\al(x)\be(x)}{\tau(x)-\al(x)}.\lab{g(delta)}
\eeq
With these notations, our task is to establish that
\beq\lab{estimateP}
N(F;\ma{P}) \ll P_3^\ve \exp\lb \frac{1}{d}g(0)\log P_3\rb,
\eeq
whenever $T \geq (P_1P_2)^d$.

Our first step is to note that $N(F;\ma{P}) \leq N(F;\ma{B})$,
and we proceed to estimate the latter.
Fortunately we shall only need to make minor alterations to the second
author's proof of
(\ref{annals}) to do so.
Let $\mcal{P} \geq \log^2 (\|F\|B_3)$, where $\|F\|$ denotes the maximum
modulus of the
coefficients of $F$.  Then according to [1; Lemma 4]
it suffices to consider the set of $\x$ counted by $N(F;\ma{B})$
for which $p \nmid \nabla F(\ma{x})$, for a fixed prime $p$
in the range $\mcal{P}\ll p\ll \mcal{P}$.
For each non-singular
$\ma{t}=(t_1,t_2,t_3)$
on the projective variety $F(\ma{t})\equiv 0~(\bmod{~p})$,
we write $S(\ma{t})$ for the set of points counted by
$N(F;\ma{B})$ for which $\x\equiv \lambda \ma{t}~(\bmod{~p})$ for some
integer $\lambda$.
Clearly there are $O(\mcal{P})$ possible values of $\ma{t}$.
Let $\delta$ be given by (\ref{del}).
We plan to show that whenever
\beq \lab{p}
\mcal{P}\gg
B_3^{\ve}
\exp\lb \frac{1}{d}g(\delta)\log B_3\rb \log^2\|F\|,
\eeq
there is an auxiliary form $G(\x)$ of degree $O(1)$,
such that  $F \nmid G$ and $G(\x)=0$ for
all $\ma{x} \in S(\ma{t})$.
It turns out that we may only do this if
\beq\lab{lowt}
\log T' \geq  d\log B_1+d\log B_2.
\eeq
Under this assumption we easily deduce the estimate
\beq\lab{estimateB}
N(F;\ma{B})\ll
B_3^{\ve}
\exp\lb \frac{1}{d}g(\delta)\log B_3\rb,
\eeq
via an application of B\'ezout's Theorem and [1; Theorem 4], just
as in the proof of (\ref{annals}).  We
now show how Theorem 1 can be derived from
(\ref{estimateB}).
We first observe that (\ref{lowt}) holds in view of the first of the
inequalities (\ref{tau}).
Next we show how (\ref{estimateP}) follows from the corresponding
estimate (\ref{estimateB}) for $N(F;\ma{B})$.
This will require the following result.

\begin{lem}\lab{1}
Let $\delta$ be given by (\ref{del}).  Then we have
$$
g(\delta) \leq g(0) +\ve.
$$
\end{lem}

Since $F$ is absolutely irreducible we must have $T \geq P_3$, and
hence (\ref{abt1}) implies the lower bound $\tau(x) \geq 1/d$.
Once combined with (\ref{tau}), we deduce that
\beq\lab{t-a}
\tau(x)-\al(x) \geq \max \{\be(x), \frac{1}{d}-\al(x)\} \geq \max
\{\al(x),
\frac{1}{d}-\al(x)\} \geq \frac{1}{2d}.
\eeq
To prove Lemma \ref{1} we shall also use the fact that for any $x \geq 0$
we have the trivial inequalities
\beq\lab{triv}
0< \al(x), \be(x), \tau(x) \leq 1, \qquad |\al'(x)|, |\be'(x)|, |\tau'(x)|
\leq 6d.
\eeq
To see the last three inequalities we note that
$$
|\tau'(x)|= \left| \frac{3}{1+\kappa x} -
\frac{\kappa(\tau+3x)}{(1+\kappa x)^2} \right| \leq 3+\kappa \leq 6d,
$$
for example, since $3 \leq \kappa \leq 3d$.
If we write
$$
h_1(x)=\frac{\al(x)+ \be(x) -\al(x)\be(x)}{\tau(x)-\al(x)},
$$
then there exists some $0 < \xi < \delta$ such that
$
h_1(\delta)-h_1(0)=\delta h_1'(\xi),
$
by the mean value theorem.
Using (\ref{tau}) it is easy to see that
$0 \leq \al(x)+\be(x)-\al(x)\be(x) \leq 1$, and so
(\ref{t-a}) and (\ref{triv}) yield
$$
|h_1'(x)| \leq  24d^2+48d^3 \leq 72d^3,
$$
for any $x \geq 0$.
Similarly, we write
$$
h_2(x)=\frac{1-\be(x)}{\tau(x)-\al(x)},
$$
and deduce that
$|h_2(\delta)-h_2(0)| \leq 60d^3 \delta$.
Finally, we write
$$
h_3(x)=\frac{\al(x)^2}{\tau(x)-\be(x)},
$$
and consider the derivative
$$
h_3'(x)=\frac{2 \al(x)\al'(x)}{\tau(x)-\be(x)} -
\lb\frac{\al(x)}{\tau(x)-\be(x)}\rb^2(\tau'(x)-\be'(x)).
$$
But (\ref{tau}) implies that $\al(x) \leq \tau(x)-\be(x)$,
and so we easily obtain  the bound $|h_3(\delta)-h_3(0)| \leq 24d\delta$.
Taken together with (\ref{g(delta)}), it therefore follows that
\begin{eqnarray*}
|g(\delta)-g(0)| &\leq&
|h_3(\delta)-h_3(0)|h_2(\delta)+
|h_2(\delta)-h_2(0)| h_3(0) + |h_1(\delta)- h_1(0)|\\
&\leq&
2d|h_3(\delta)-h_3(0)|+
|h_2(\delta)-h_2(0)| + |h_1(\delta)- h_1(0)|\\
&\leq&
\delta \{48d^2+ 60d^3 + 72d^3\}\\
&\leq&
\ve,
\end{eqnarray*}
by (\ref{del}).  This completes the proof of Lemma \ref{1}.

We are now in a position to deduce (\ref{estimateP}) from
(\ref{estimateB}) and Lemma \ref{1}.
Using the same notation  it is easy to deduce from
(\ref{abt1}), (\ref{tau}) and (\ref{g(delta)}) that
$$
g(0)= -h_3(0)h_2(0) + h_1(0) \leq 4d.
$$
Recall that  $\kappa \leq 3d$.
Then (\ref{b}), (\ref{del}) and Lemma \ref{1} yield
\begin{eqnarray*}
\frac{1}{d}g(\delta)\log B_3 - \frac{1}{d}g(0)\log P_3 &\leq&
\frac{1}{d}\log P_3\{\ve+\kappa \delta(g(0)+\ve)\}\\
&\leq& 2\ve \log P_3,
\end{eqnarray*}
whence (\ref{estimateP}).

We now show how the lower bound (\ref{p}) suffices for the existence
of a suitable auxiliary form $G(\x)$.
Let $D \geq d$ and $A>0$, and consider the exponent set
$$
\mcal{E}(A) = \left\{\ma{e}\in \Z^3: e_i \geq 0, ~\sum_{i=1}^3 e_i=D,
~\sum_{i=1}^3
  e_i\log B_i \leq A, ~\exists j \rom{~s.t.~}
e_j< f_j\right\}.
$$
Let $E=\#\mcal{E}(A)$, and suppose for the moment that $E\leq
\#S(\ma{t})$.  If we choose any distinct vectors
$\x^{(1)},\ldots,\x^{(E)} \in S(\ma{t})$, then
it will actually suffice just to show that the determinant
$$
\D=\det(\x^{(i)\ma{e}})_{1\leq i\leq E,\;\ma{e}\in\mcal{E}(A)}
$$
vanishes whenever (\ref{p}) occurs.   Indeed, the construction of the
auxiliary polynomial
$$
G(\x)=\sum_{\ma{e}\in \mcal{E}(A)}
a_\ma{e} x_1^{e_1}x_2^{e_2}x_3^{e_3}
$$
is then identical to that given in the proof of (\ref{annals}).
We have written
$
\ma{w}^{\ma{e}}=w_{1}^{e_{1}}w_2^{e_2}w_{3}^{e_3}
$
in the definition of $\D$, in which rows correspond to the different
vectors $\x^{(i)}$,  and columns
correspond to the various $\ma{e} \in \mcal{E}(A)$.
Furthermore,  it is immediate from
the proof of (\ref{annals}) that any such form
$G$  cannot be divisible by $F$.
The advantage over the previous situation is that our new exponent
set $\mcal{E}(A)$ allows us to choose an optimal value of $A$, for
which better control over the size of $|\D|$ is possible
in certain situations.

Our proof now breaks into two parts.  Firstly we must obtain an
estimate for the real modulus of  $\D$, and then
secondly show that its $p$-adic order is sufficiently large that $\D$
must in fact vanish.
We begin with the first of these, and use the fact that
$|x_j^{(i)}| \leq B_j$ for $1\leq j \leq 3$ to deduce that the column
corresponding to the exponent vector $\ma{e}$ consists of elements
of modulus at most
$
B_1^{e_1}B_2^{e_2}B_3^{e_3}.
$
It therefore follows that
\beq\lab{delta}
|\Delta|\leq E^{E}\prod_{\ma{e}\in\mcal{E}(A)}
B_1^{e_1}B_2^{e_2}B_3^{e_3}.
\eeq
For any $\ma{e}\in\Z^3$ with $e_i\ge 0$, we henceforth set
$$
\sigma(\ma{e})=e_1\log B_1+e_2\log B_2+e_3\log B_3.
$$
For $1\leq i \leq 3$, define $\mcal{E}_i$ to be the subset of
$\mcal{E}(A)$ for which
$e_i<f_i$.   Then we have
\beq\lab{start}
\log \prod_{\ma{e}\in\mcal{E}(A)} B_1^{e_1}B_2^{e_2}B_3^{e_3}
=\sum_{\ma{e}\in\mcal{E}(A)}\sigma(\ma{e})\leq \sum_{i=1}^3
\sum_{\ma{e}\in\mcal{E}_i}\sigma(\ma{e}),
\eeq
since $\mcal{E}(A)=\mcal{E}_1\cup \mcal{E}_2\cup
\mcal{E}_3$.   Hence it suffices to estimate
$\sum_{\ma{e}\in\mcal{E}_i}\sigma(\ma{e})$, for which it will be
convenient to write
$$
b_i=\log B_i, \qquad 1 \leq i \leq 3.
$$
Let $c=1/\delta+\kappa$,
where $\kappa$ and $\delta$ are given by (\ref{kappa}) and
(\ref{del}), respectively.
Then it follows from our initial change of variables (\ref{b})
that we have the inequalities
\beq\lab{ass}
0< b_1 < b_2 <  b_3 \leq c b_1, \qquad b_3 \leq c(b_j-b_i),
\eeq
for each $1 \leq i<j \leq 3$.

It is convenient at this point to make the assumption that
$A$ is contained in
the interval
\beq\lab{Arange}
D b_2 < A \leq D b_3.
\eeq
Clearly our definition of $\mcal{E}(A)$ would be rather pointless
if we allowed $A > Db_3$, since then the condition
$\sigma(\ma{e}) \leq A$ is automatic and we retrieve
the exponent set considered in the proof of (\ref{annals}).
Similarly,
$\mcal{E}(A)$ is obviously empty for any $A \leq Db_1$.  Our motive
for omitting any treatment of the interval $Db_1 < A \leq Db_2$
is not so apparent.  Indeed, it is possible to adjust our
argument to take such values of $A$ into account and actually
achieve something new at the end of it.  We have chosen not to do so
simply because we would
ultimately be led to a weaker result than Theorem 1.  Moreover,
we are able to simplify our work considerably under the
hypothesis  (\ref{Arange}).

Henceforth let $i,j,k$ denote distinct elements of the set
$\{1,2,3\}$, and define
$$
\mcal{M}_{jk} = \left\{(m_j,m_k)\in \Z^2: m_j,m_k \geq 0, ~m_j+m_k=D,
~\tau(m_j,m_k) \leq A\right\},
$$
where $\tau(m_j,m_k)=m_jb_j+m_kb_k$.
We shall apply the following result to simplify our estimate for
$\sum_{\ma{e}\in\mcal{E}_i}\sigma(\ma{e})$.

\begin{lem}\lab{2}
For $1 \leq i \leq 3$ we have
$$
\sum_{\ma{e}\in\mcal{E}_i}\sigma(\ma{e}) \leq
f_i\sum_{(m_j,m_k)\in\mcal{M}_{jk}}\tau(m_j,m_k)  +O(A),
$$
provided that (\ref{Arange}) holds.
\end{lem}

We prove the result for $i=3$, say, and consider  values of $\ma{e} \in
\mcal{E}_3$ with a fixed component $e_3<f_3$.
If $m_1=e_1$ and $m_2=e_2+e_3$ then
$|\tau(m_1,m_2)-\sigma(\ma{e})| \leq db_3$. In particular we either
have $(m_1,m_2) \in \mcal{M}_{12}$ or
\beq\lab{r12}
A<\tau(m_1,m_2) \leq A+d b_3.
\eeq
Let $R_{12}$ denote the number of non-negative $m_1,m_2\in \Z$ such that
$m_1+m_2=D$ and (\ref{r12}) holds.
It follows that
$$
\sum_{\ma{e}\in\mcal{E}_3}\sigma(\ma{e}) \leq
\sum_{e_3<f_3}\left\{A R_{12}+ 
\sum_{(m_1,m_2)\in\mcal{M}_{12}}(\tau(m_1,m_2)+db_3)
\right\}.
$$
Moreover, we observe that
\beq\lab{interval}
A>D b_2 \gg D b_3,
\eeq
by (\ref{ass}) and (\ref{Arange}).
Since $\#\mcal{M}_{12} \leq D+1$ and $D \geq d$, this yields
$$
\sum_{\ma{e}\in\mcal{E}_3}\sigma(\ma{e}) \leq
f_3 \sum_{(m_1,m_2)\in\mcal{M}_{12}}\tau(m_1,m_2)+ O(A+AR_{12}).
$$
In order to estimate $R_{12}$ we set $\nu=b_1D/(b_2-b_1)$.  Then
$R_{12}$ is the number of non-negative integers $m_2 \leq D$ for which
$$
\frac{A}{b_2-b_1}< m_2+\nu \leq \frac{A}{b_2-b_1}+\frac{db_3}{b_2-b_1}.
$$
It follows from (\ref{ass}) that $db_3/(b_2-b_1) \ll 1$, whence
$R_{12}=O(1)$.   This  suffices for the proof of
Lemma \ref{2}.

Recall that $\{i,j,k\}$ is a permutation of $\{1,2,3\}$.
Our next task is to establish the estimate
\beq\lab{claim1}
\sum_{(m_j,m_k)\in \mcal{M}_{jk}} \tau(m_j,m_k) =
\frac{1}{2}\sum_{g,h}\mbox{}^{*}\frac{A^2-D^2b_h^2}{b_g-b_h} +O(A),
\eeq
where $\{g,h\}$ runs over permutations of $\{j,k\}$, and
$\Sigma^*$ denotes the condition $b_h < A/D$.
Taking the case corresponding to $i=3$ first,
we automatically have $\tau(m_1,m_2)<A$ since $Db_2<A$, by
(\ref{Arange}).   Therefore
\begin{eqnarray*}
\sum_{(m_1,m_2)\in \mcal{M}_{12}} \tau(m_1,m_2) &=&
\sum_{m_1+m_2=D}\tau(m_1,m_2)\\
&=& \sum_{m_2=0}^D \{(b_2-b_1)m_2 +b_1D\}\\
&=& (b_1+b_2)D(D+1)/2\\
&=&
\frac{1}{2}\left\{\frac{A^2-D^2b_1^2}{b_2-b_1}+\frac{A^2-D^2b_2^2}{b_1-b_2}\right\}
+O(A).
\end{eqnarray*}
In the case $i=1$, we find
\begin{eqnarray*}
\sum_{(m_2,m_3)\in \mcal{M}_{23}} \tau(m_2,m_3) &=&
  \sum_{m_3=0}^{[\frac{A-Db_2}{b_3-b_2}]} \{(b_3-b_2)m_3 +b_2D\}\\
&=& \frac{A^2-D^2b_2^2}{2(b_3-b_2)} +O(A),
\end{eqnarray*}
on using (\ref{interval}).  The treatment of the case $i=2$ is similar.
This completes the proof of (\ref{claim1}).

Upon combining (\ref{delta}), (\ref{start}), Lemma \ref{2} and
(\ref{claim1}), we are therefore
led to the following result.

\begin{lem}\lab{estimateD}
We have
$$
\log |\Delta| \leq
\frac{1}{2}\sum_{i=1}^3 f_i
\sum_{g,h}\mbox{}^{*}\frac{A^2-D^2b_h^2}{b_g-b_h} + E\log E+ O(A),
$$
provided that (\ref{Arange}) holds.
\end{lem}

For our fixed prime $p$ with order of magnitude $\mcal{P}$, we must
now determine a lower bound for $\nu_p(\D)$, the $p$-adic order of $\D$.
However it suffices to compute $E=\#\mcal{E}(A)$,
since [1; Lemma 6] implies that
\beq\lab{nu}
\nu_p(\D) \geq
\frac{1}{2}E^2 \{1+o(1)\},
\eeq
as $E \rightarrow \infty$.
In fact a lower bound for $E$ can easily be deduced by mimicking the
previous
calculation.  Beginning
with the analogue of (\ref{start}), one obviously has
\beq\lab{start-lower}
\left|
\sum_{\ma{e}\in\mcal{E}(A)}1 - \sum_{i=1}^3
\sum_{\ma{e}\in\mcal{E}_i} 1 \right|
\leq
\sum_{i<j}\sum_{\ma{e}\in\mcal{E}_i\cap \mcal{E}_j}1 =O(1).
\eeq
Moreover, the corresponding version of  Lemma \ref{2} is
\beq\lab{lemma1-lower}
\sum_{\ma{e}\in\mcal{E}_i}1 \geq
f_i\sum_{(m_j,m_k)\in\mcal{M}_{jk}}1+O(1),
\eeq
under the same assumption that (\ref{Arange}) holds.  We prove this for
$i=3$, say.
Recall that
$m_1+m_2=D$ and $\tau(m_1,m_2) \leq A$ whenever $(m_1,m_2) \in
\mcal{M}_{12}$. In particular, for each integer $0 \leq k <f_3$ we
either have $(m_1,m_2-k,k)
\in \mcal{E}_3$ or  $0 \leq m_2 <k$ or
$A-db_3< \tau(m_1,m_2) \leq A$.
Hence an argument similar to that used in the proof of Lemma \ref{2}
yields the upper bound
$$
\sum_{k<f_3}\sum_{(m_1,m_2)\in\mcal{M}_{12}} 1 \leq
\sum_{\ma{e}\in\mcal{E}_3}1 +O(1),
$$
which establishes (\ref{lemma1-lower}).
Combining (\ref{start-lower}) and (\ref{lemma1-lower}), together with
the corresponding version of (\ref{claim1}), we obtain the
lower bound
\beq\lab{lowerE}
E \geq \sum_{i=1}^3 f_i
\sum_{g,h}\mbox{}^{*}\frac{A-Db_h}{b_g-b_h} +O(1),
\eeq
provided that (\ref{Arange}) holds.

Let $A$ be contained in the interval (\ref{Arange}).  Then as $D
\rightarrow \infty$ we will have both $E \rightarrow \infty$ and
$A=o(D^2b_2)$.  It now follows from Lemma \ref{estimateD} that
$$
\log |\Delta| \leq
\frac{1}{2}\left\{f_1 \frac{A^2-D^2b_2^2}{b_3-b_2}+f_2
  \frac{A^2-D^2b_1^2}{b_3-b_1}+f_3 D^2(b_1+b_2)\right\}(1+o(1)) +o(E^2),
$$
as $D \rightarrow \infty$.  Moreover from (\ref{nu}) and
(\ref{lowerE}) we also have
$$
\nu_p(\D) \geq
\frac{1}{2}
\left\{f_1 \frac{A-Db_2}{b_3-b_2}+f_2
  \frac{A-Db_1}{b_3-b_1}+f_3 D\right\}^2(1+o(1)).
$$
We may therefore conclude that
$$
\frac{\log |\Delta|}{\nu_p(\Delta)} \leq \frac{f_1
\frac{A^2-D^2b_2^2}{b_3-b_2}+f_2
  \frac{A^2-D^2b_1^2}{b_3-b_1}+f_3 D^2(b_1+b_2)}{\left\{f_1
\frac{A-Db_2}{b_3-b_2}+f_2
  \frac{A-Db_1}{b_3-b_1}+f_3 D\right\}^2}(1+o(1)),
$$
as $D \rightarrow \infty$.
Define the constants
$$
\lambda
 =\frac{db_3-\log T'}{(b_3-b_1)(b_3-b_2)}, \qquad
  \phi
  = \frac{db_1b_2+b_3(\log T'-db_1-db_2)}{(b_3-b_1)(b_3-b_2)},
$$
and
$$
\gamma
= \phi(b_1+b_2)+\lambda b_1b_2.
$$
Then in particular $\lambda \geq 0$ and (\ref{lowt}) implies that
$\phi > 0$.
We shall consider the behaviour of the real-valued function
$$
f(A)= \frac{\lambda A^2+\gamma D^2}{(\lambda A+\phi D)^2},
$$
as $A$ varies over the interval (\ref{Arange}).
We recall that $\D$ necessarily vanishes if
$p^{\nu_p(\D)} > |\D|$.  Using the identities $f_1+f_2+f_3=d$ and
$\sigma(\ma{f})=\log
T'$, a straightforward calculation reveals that $\D$ vanishes if
$\log p>(1+o(1))f(A)$ for any $A$ in the interval (\ref{Arange}).

It remains to choose a suitable value of $A$ for which
the function $f(A)$ is minimized.
In fact, an easy calculation reveals that $f$ has a turning point at
$A=\gamma D/\phi$.  Moreover,
this value of $A$ is contained in  the interval (\ref{Arange}) precisely
when
$\phi b_2< \gamma \leq \phi b_3$.    The lower bound here always
holds, whereas it is not hard to see that the upper bound is true
if and only if (\ref{lowt}) holds.
Hence it suffices to take
$$
p >
\exp\lb \frac{db_1b_2b_3 +(b_1b_3+b_2b_3-b_1b_2)(\log
  T'-db_1-db_2)}{(\log T' - db_1)(\log T' - db_2)}(1+o(1))\rb,
$$
under this assumption.  Therefore (\ref{p}) is indeed satisfactory,
provided that
$D$ is chosen to be sufficiently large in terms of $\ve$ and $d$.

\section{Equal sums of two powers}

Let $k \geq 4$ and $X \geq 1$.
We now turn to our estimate for the number $N_{k}^{(0)}(X)$ of positive
non-trivial integral
solutions of the Diophantine equation (\ref{x}), which are contained
in the region  $\max \{w,x,y,z\}\leq X$.
It will suffice to count positive integers $w,x,y,z$ such that
$x <y \leq z<w$.
For each such solution we define
$$
v_1=z-x, \qquad v_2=z+x.
$$
It follows that
\beq\lab{region}
1\leq y < w \leq X, \qquad 1 \leq v_1<v_2 \leq 2X.
\eeq
Furthermore, under this transformation (\ref{x}) takes the shape
\beq\lab{uv}
2^{k-1}\{w^k-y^k\}=v_1 f(v_1,v_2),
\eeq
where
\beq\lab{f}
f(v_1,v_2)=
\sum_{0 \leq j <k/2} \left( k \atop 2j+1 \right) v_1^{2j}v_2^{k-2j-1}
\eeq
is a binary form of degree $k-1$.
In order to estimate the number of integers $w,y,v_1, v_2$ such that
(\ref{region}) and (\ref{uv}) hold, we begin by considering the
contribution corresponding to a fixed choice of
$v_1$.  For this we define
$$
\xi=\xi(v_1) =\prod_{p \mid v_1, p>2} p
$$
to be the odd square-free kernel of $v_1$, and consider the set
$$
S= \{(w,y) \in \Z^2: \xi \mid w^k-y^k  \}.
$$
Let $(w,y) \in S$ and let $p$ be any prime divisor of $\xi$. Then we
see that either $p$ divides $y$, or there exist at most $k$ 
integers
$\lambda_1, \ldots, \lambda_t$, say,  such that
$$
w \equiv \lambda_i y \pmod{p}
$$
for some $1 \leq i \leq t$.
Collecting these lattice conditions together
via the Chinese Remainder Theorem, we therefore conclude
that $S$ is a union of $O(k^{\omega(\xi)})$ lattices in $\Z^2$,
each of determinant $\xi$.

We henceforth fix our attention upon those
$w,y$ contained in the region (\ref{region}), which lie
in one such lattice $\sfl$, say.
By \cite[Lemma 1, (iii)]{annal} there exist basis vectors
$\ma{e}^{(1)}, \ma{e}^{(2)} \in \sfl$ with
\beq\lab{det}
\xi \ll |\ma{e}^{(1)}||\ma{e}^{(2)}| \ll \xi,
\eeq
and such that whenever we write
\beq\lab{subs}
(w,y)=u_1\ma{e}^{(1)}+u_2\ma{e}^{(2)}
\eeq
for appropriate integers $u_1, u_2$, we automatically have
\beq\lab{region'}
u_1 \ll X/|\ma{e}^{(1)}|, \qquad u_2 \ll X/|\ma{e}^{(2)}|.
\eeq
We assume without loss of generality that $|\ma{e}^{(1)}| \leq
|\ma{e}^{(2)}|$,  and so
\beq\lab{e}
|\ma{e}^{(1)}|^i|\ma{e}^{(2)}|^j \geq
(|\ma{e}^{(1)}||\ma{e}^{(2)}|)^{(i+j)/2}
\eeq
for any choice of $j \geq i \geq 0$.

For each $v_1$, it
suffices to count the number of integers $u_1,u_2, v_2$ lying in the
region defined by
(\ref{region}) and (\ref{region'}), for which
$$
g(u_1,u_2)
=v_1f(v_1,v_2),
$$
where
$g(u_1,u_2)$ is obtained from $2^{k-1}\{w^k-y^k\}$ via the
substitution (\ref{subs}).
In fact we shall also fix a choice of $u_2$, and then count the number
$M(X;u_2,v_1)$, say,
of integers $r, s$ for which
$r \ll X/|\ma{e}^{(1)}|$, $1 \leq s \leq 2X$, and
\beq\lab{pq}
p(r)=q(s),
\eeq
where
$p(r)=g(r,u_2)$ and $q(s)=v_1f(v_1,s)$.

An important issue here is whether or not the  polynomial
$p-q$ is absolutely irreducible.
But it is not hard to see that (\ref{pq}) is obtained from (\ref{x})
via an appropriate affine plane section.  We now distinguish the
projective plane sections of
(\ref{x}) into three distinct types: those producing absolutely 
irreducible curves of degree $k$, those
that produce a line and an absolutely irreducible curve
of degree $k-1$, and finally those that produce a union of absolutely
irreducible curves each of degree $\leq k-2$.
In the first case it is clear that the corresponding affine plane
section is absolutely irreducible.  In the second case we
may deduce that the polynomial defining the resulting affine curve
(\ref{pq}) is either the
product of a linear polynomial and an absolutely irreducible
polynomial of degree $k-1$, or it is absolutely irreducible of
degree $k-1$.
The latter possibility is clearly satisfactory, whereas the former
possibility implies that $p-q$ has degree $k$.  But then
$\deg p=k$ and $\deg q=k-1$, and 
a result of Schmidt \cite[Theorem III.1B]{schmidt} tells us that 
$p(r)-q(s)$ should be absolutely irreducible.

In the final degenerate case, we may assume that
the projective plane section of (\ref{pq}) produces
at least two distinct absolutely
irreducible curves.
Indeed, the existence of any plane section producing
precisely
one line in (\ref{x}) would imply the existence of a singular point
on the surface.
Now we know by the
previously discussed result of Colliot-Th\'el\`ene \cite[Appendix]{annal}
that
(\ref{x}) contains finitely many plane curves of degree $\leq k-2$.
It follows that there can only be $O(1)$ projective plane sections
which produce two distinct absolutely irreducible curves of degree $\leq
k-2$.  Recall from the introduction that any absolutely irreducible plane
curve of degree $d$
contains $O(X^{2/d+\ve})$ rational points of height at most
$X$.  Since trivial integral solutions to (\ref{x})
correspond to rational points lying on projective lines in the surface,
we therefore conclude that
there is a total contribution of $O(X^{1+\ve})$ to
$N_k^{(0)}(X)$, from those affine plane sections of (\ref{x}) which
lead to reducible curves (\ref{pq}).
This is clearly satisfactory for Theorem $2$, and we may assume henceforth
that $p-q$ is
absolutely irreducible.

In order to estimate $M(X;u_2,v_1)$, we shall apply
(\ref{cime}) to
$p(r)-q(s)$.
In view of the shape (\ref{f}) that $f$ takes, it is apparent
that $q(s)$ contains the monomial $s^{k-1}$ with non-zero coefficient.
Taking $P_2\ll X/|\ma{e}^{(1)}|$ and $P_3=2X$ in (\ref{cime}), we see that
$T \gg X^{k-1}$ and hence
$$
M(X;u_2,v_1) \ll X^{1/(k-1)+\ve}|\ma{e}^{(1)}|^{-1/(k-1)}.
$$
Summing over the values of
$u_2 \ll X/|\ma{e}^{(2)}|$ in (\ref{region'}) we therefore obtain the
contribution
\begin{eqnarray}
&\ll&
\frac{X^{1+1/(k-1)+\ve}}{|\ma{e}^{(1)}|^{1/(k-1)}|\ma{e}^{(2)}|}\nonumber\\
&\ll&
(\xi^{-1/2}X)^{1+1/(k-1)+\ve},\lab{tmp1}
\end{eqnarray}
via (\ref{e}) and (\ref{det}).
Let $Y \geq 1$, and write
$A_\ve=(1-2^{-\ve})^{-1}$ for a fixed choice of $\ve>0$. Then for
any $\theta\leq 1$ we have
\begin{eqnarray*}
\sum_{n \leq Y} \xi(n)^{-\theta}
&\leq& \sum_{n \leq Y} \xi(n)^{-\theta} \lb
\frac{Y}{n}\rb^{1-\theta+\ve}\\
&=& Y^{1-\theta+\ve} \sum_{n \leq Y} \xi(n)^{-\theta} n^{\theta-1-\ve}\\
&\leq& Y^{1-\theta+\ve} \sum_{n=1}^\infty \xi(n)^{-\theta}
n^{\theta-1-\ve}\\
&\leq& Y^{1-\theta+\ve} \prod_{p} \left\{1+p^{-\theta}
  p^{\theta-1-\ve}+p^{-\theta} p^{2\theta-2-2\ve}+\ldots\right\}\\
&=& Y^{1-\theta+\ve} \prod_{p} \left\{1+p^{-\theta}
  \frac{p^{\theta-1-\ve}}{1-p^{\theta-1-\ve}}\right\}\\
&\leq& Y^{1-\theta+\ve} \prod_{p} \left\{1+p^{-1-\ve}A_\ve\right\}\\
&=& c_\ve Y^{1-\theta+\ve},
\end{eqnarray*}
say.
Upon taking $\theta= \frac{1}{2}(1+1/(k-1))$, so that $\theta \leq 1$ for
every $k \geq 2$, it
therefore follows from (\ref{tmp1}) that
\begin{eqnarray*}
N_{k}^{(0)}(X) &\ll& X^{1+1/(k-1)+\ve} \sum_{v_1 <2X} \xi(v_1)^{-\theta}\\
&\ll  & X^{3/2+1/(2k-2)+\ve}.
\end{eqnarray*}
This completes the proof of Theorem 2.

\end{document}